\documentclass[12pt,a4paper]{article}
\usepackage{t1enc}
\usepackage[latin1]{inputenc}
\usepackage[english]{babel}
\begin{document}

\date{}
\title{Ricci Curvature  and Singularities of Constant Scalar Curvature Metrics}
\vspace{1cm}

\author{A. Raouf  Chouikha\footnote{Chouikha@math.univ-paris13.fr}
}

\maketitle

\vspace{1cm}
\begin{center}
 University of Paris-Nord \\ LAGA,CNRS UMR 7539 \\ Av. J.B. Clement \\
 Villetaneuse F-93430 

\vspace{7cm}

\end{center}
\newpage
\begin{abstract}
In this work we consider periodic spherically symmetric metrics of constant positive scalar curvature on the n-dimensional cylinder called pseudo-cylindric metrics. These metrics belong to the conformal class $[g_0]$ of the Riemannian product $S^1\times S^{n-1}$ : a circle of length $T$ crossed with the (n-1)-dimensional standard sphere.\\
Such metrics have a harmonic Riemannian curvature and a non parallel Ricci tensor, except for the cylindric one. Thus, it appears a natural link between them and the Derdzinski metrics which are warped product and classify a family of Riemannian manifolds. These two families actually differ by conformal transformations. Moreover, we are interested in the multiplicity problem of the pseudo-cylindric metrics in $[g_0]$. We also study the existence problem and the number of the Derdzinski metrics. \\
Furthermore, we prove that the pseudo-cylindric metrics may be expressed in terms of elliptic functions for the dimension $n = 3,4$ and $6$ only, and in terms of automorphic functions for any other dimension. This fact allows us to give new bounds for global estimates. \\
Finally, we examine the curvature of the asymptotic pseudo-cylindric metrics which are complete singular Yamabe metrics on the standard sphere punctured of $k$ points. We show that any of such (non trivial) metric has a non parallel Ricci tensor.

\end{abstract}
\vspace{5cm}

{\it Mathematical Subjects Classification : 53C21, 53C25, 58G30.}

\newpage

\tableofcontents

\newpage
\section {Introduction and statements of results }
\smallskip

  Let us consider a smooth compact Riemannian manifold\ $(M,g)$\ of dimension \
$n\geq 3$ , where \ $g$ \ is a smooth Riemannian metric on \ $M.$  \ $dv$\ denotes
the volume form of \ $g$,\ and \ $R(x)$ \ its  scalar curvature. \    
 $ {\mathcal H_1}^2(M,g)$\ is the Sobolev space.\\
For any  non null function\quad $u \in {{\mathcal H_1}^2(M,g)},$\quad  we
define the functional 
\begin{equation}
 \qquad J(u) = {{{4{{n-1}\over{n-2}}}\int_M\left|\nabla
u\right|^2 dv + \int_M R(x) u^2 dv}\over{\left(\int_M u^{2n\over{n-2}}dv\right
)^{n-2\over n}}}\  .
\end{equation}
The critical points of the functional \ $J(u)$ \ are solutions of the
Euler-Lagrange equation 
\begin{equation}
 \qquad 4{{n-1}\over{n-2}} \Delta{_g} u + R(x) u = k
|u|^{4\over{n-2}} u\ ,  
\end{equation}
where \ $k$ \ is a real number.\\

There is a geometric interest only if we assume the solution \ $u$ \ of (2) to be \ $C^\infty$ \ positive. In
this case, the conformal metric \ $g' = u^{4\over{n-2}} g$ \ has  constant
scalar curvature  \ $k$ \ if \ $u$ \ is a critical point of the
functional \ $J(u)$\ and conversely. So, we know by resolution of the Yamabe 
conjecture for the smooth compact manifolds, the infimum of\ $J(u)$\ is always attained . A such  metric with constant scalar curvature belonging to the conformal class \ $[g]$\ is called Yamabe metric. \\ More precisely, let \ $(S^n,g_0)$ \ be the standard sphere, we may assert, [L-P] : \\ 

{\it There exists a  non zero function \ $u_0$, \ positive, \ $C^\infty$ on\
$(M,g)$\ such that  $$J(u_0) = Inf_{u\in H_1^2(M)} J(u) = \mu (M,g) .$$
\quad  Moreover, we have  \ $\mu (M,g) < \mu (S^n,g_0)$,\  equality holds
 only if there exists a conformal diffeomorphism between the manifold \
$(M,g)$\ and the standard sphere\ $(S^n,g_0)$.}\\

\smallskip
 \quad It is also known that the Yamabe metric is unique up to homothety, when
\ $(M,g)$\ is such that \quad $\int_M R(x) dV \leq 0$ \ or when \ $(M,g)$ is an
Einstein Manifold. \\ For the positive case of the scalar curvature \ ($k > 0$),\ it is possible to get multiple Yamabe solutions for certain conformally flat product Riemannian manifolds. \\

 The multiplicity problem of the Yamabe
solutions have been intensively studied during the last few years,  in the more
general context known as the singular Yamabe problem with isolated singularities. The latter consists to find a
metric \ $g = u^{4\over{n-2}} g_0$ \ on the domain \ $S^n - {\Lambda_k}$,\ where \
${\Lambda_k}$\ is a finite point-set of \ $(S^n,g_0).$\ This metric being complete and have a
constant scalar curvature
\ $ R(g).$

\smallskip

 Concerning the special case \ $k = 2$,\  R. Schoen (see section 2 of\  [Sc]) has partially described the conformal class
of the product metric
\
$(S^1
\times S^{n-1},dt^2+d{\xi}^2), $ \ where \ $S^1$\ is the circle of length\ $T$\ and\ $(S^n,g_0)$\ 
is the standard sphere.

  He proved the existence and the number of their rotationally invariant Yamabe metrics\quad
$${g^j}_c = {{u^j}_c}^{4\over{n-2}}(dt^2+d{\xi}^2)\ ;\quad j = 1,2,...,k-1 $$ called pseudo-cylindric
metrics (sometime also called by different authors Delaunay metrics or Fowler metrics) which depend on T-values. But, his proof is not completed. The completion has not been yet published elsewhere ( unless our unpublished joint work with F. Weissler [Ch-W]).\\
Now, it is well known that the number of these metrics  is finite in a conformal class for the product \ $[dt^2+d{\xi}^2]$.\\ 
We refer to the papers [Sc], [Ch3] or [Ch-W].\\
 In particular, we get a unique Yamabe metric if 
\begin{equation}
\qquad T \leq  2\pi\sqrt{{n-1}\over{R_0}}\ .
\end{equation}
This metric is \ $c^{te}(dt^2+d{\xi}^2)$.\\

 One of our main results concerns geometric properties of the above metrics. That is any non trivial pseudo-cylindric metric
 have a Ricci tensor non parallel and so, its Riemannian curvature is harmonic (see
Theorem (4-3)below and [Ch3]). A tensor of order two \ $\tau = (\tau_{ij})$\ is said to be nonparallel if the
covariant differentiation of\ $\tau$\ along \ $(M,g)$\ has a nonzero component.\\

 Observe that all these metrics -sharing the above curvature properties- have been
classified by A. Derdzinski. Indeed, these metrics are covered  by a warped product metrics :
\quad
$ dt^2+{f^2}(t)d{\xi}^2$  \ up to an isometry, where  \ $f(t) $\quad satisfy a certain ODE. Moreover, we show  by using a
Yau result (Lemma (4-5) below), that any pseudo-cylindric metric may actually be identified 
 with a Derdzinski metric, up to a conformal
diffeomorphism. This identification may be written
$$dt^2+{f^2}(t)d{\xi}^2 = {F^*}({{u^j}_c}^{4\over{n-2}}(dt^2+d{\xi}^2)),$$ where \ $F$\ is a
conformal diffeomorphism. \\ 
Remark that the metrics originally contructed by A. Derdzinski, permit to provide
example of compact Riemannian manifold with harmonic curvature and non-parallel Ricci tensor. We know that these  metrics are depending on the 
\ $S^{1}$- parametrization \ $T$\ (thus, there are obviously not conformally equivalent)  .\\

Nevertheless, the pseudo-cylindric metrics belong to the
same conformal class which is  parametrized by
\ $T$.\  Hence, they are conformally equivalent to the cylindric metric. That is not the same
for the Derdzinski metrics. For this reason, from the multiplicity point of view of
the Yamabe solutions, the pseudo-cylindric metrics are more interesting. This fact will be clarified below.\\ 
\\
This paper is organized as follows :\\

{\bf 1 -} We achieve in Section 2, the completed  proof (not yet  published somewhere) concerning the existence and the finite number of the pseudo-cylindric metrics in a conformal class. \\

{\bf 2 -} In section 3, we consider the link between the multiplicity of Yamabe problem 
and the best constant of sobolev. In particular, we put right and improve a theorem of T. Aubin.\\

{\bf 3 -} In Section 4, we examine the non parallelism properties of the Ricci curvature for the pseudo-cylindric metrics.\\

{\bf 4 -} In the next Section 5 we are interested on the existence of many Derdzinski metrics, satisfying certain   conditions. In particular, we show the existence
of a finite number (depending on the circle 
length  $T$ ) of such metrics. \\
This fact has not been remarked in the Part D of chapter 16 of A. Besse [Be]. His Theorem 16.38 only proved existence of non trivial Derdzinski metrics without precising the interconnection between these metrics and the circle length.\\ In particular, we exhibit a constant \
$T_0$ \ such that, if we assume\quad $T \leq  T_0$\ , then there is only one Derdzinski metric with
parallel Ricci tensor, i.e. the trivial one (see Theorem (5-3) below).\\

{\bf 5 -} Furthermore, for the dimension values 
\ $n = 3,4,6$.\ the 
pseudo-cylindric solution can be  expressed in terms of a real periodic elliptic functions. 
 The interest of the latter remark is that an elliptic function has only polar singularities. So, as a corollary (see also our previous work [Ch1]), 
 we deduce the estimates for
the corresponding pseudo-cylindric solutions (Theorem (6-1) below):
$${\displaystyle  r^{{n-2}\over 2} u(r) \leq  C,\qquad and\qquad r^{n\over 2} {{du(r)}\over {dr}}} \leq 
C'\ .$$ 
This fact allows us to explicitly determine the bounds which naturally depend  on the dimension value \ $n$ \ and the period \ $T$.\\
This generalizes a result of Cafarelli-Gidas-Spruck [C-G-S] - 
According to Pollack (Theorem (4-3) of [P] ), the first inequality  has been etablished by R. Schoen - .\\  Moreover, we show for any dimension other than \ $3,4$ \ or \ $6$, the pseudo-cylindric solutions may be expressed in terms of automorphic functions. So, they are invariant under a substitution group.\\

{\bf 6 -} Finally in section 7, we consider the asymptotic pseudo-cylindric metrics.
They are complete Yamabe metrics on \ $S^n - {\Lambda_k}$,\ where \
$\Lambda_k = \{p_1,p_2,...,p_k\}$\ (here \ $k > 2 $).\ In particular, we show that for any non trivial 
metrics, its Ricci tensor is also non parallel (Theorem (7-3) below).\\

All the proofs are exposed to the end, in Section 8.

\newpage

\section{Singular Yamabe problem and pseudo-cylindric metrics}
\subsection{The pseudo-cylindric metrics}

\quad Let a finite point-set\ $\Lambda_k = {p_1,p_2,...,p_k}$\ of the standard sphere\ $S^n$.\
$R_0$\ being its scalar curvature. consider the following\\ 

   {\bf Singular Yamabe problem:}\\ 
{\it There exists a metric \ $g = u^{4\over{n-2}} g_0$ \ in the conformal class\ $[g_0]$\ on the domain 
\  $S^n - {\Lambda_k}$,\  which is complete and
has a constant scalar curvature \ $ R(g).$} \\
In others words, the following differential equation
\begin{equation}
 4{{n-1}\over{n-2}} \Delta_{g_0} u + R_{g_0} u - R_g
u^{{n+2}\over{n-2}} = 0,    
\end{equation} 
 $$g=u^{4\over {n-2}}g_0 \quad \mbox{is complete on} \quad
S^n-\Lambda\ ,
\qquad R_g = \mbox{constant} > 0,$$  
has a solution.\\

Let \quad ${\mathcal L}_g = \Delta_g + {n-2\over{4(n-1)}} R_g$ \quad denotes 
the conformal Laplacian of any Riemannian manifold \ $(M,g).\quad {\mathcal L}_g$ \ posseses
a special property:\\

{\it Invariance property of the conformal Laplacian : Let \ $(M,g)$ \ and $(N,h)$ two Riemannian 
manifolds, such that there exists a conformal diffeomorphism\ $f : 
M \rightarrow N$ (that means if \ $f$ \ is a diffeomorphism,  there is a function
\ $v : M \rightarrow IR^+,$ \ such that \ $f^*h = v^{4\over{n-2}} g$).\\
Thus, for any function \ $\varphi : N \rightarrow IR$ \ of class \ $C^2$,\ we have
$${\mathcal L}_h(\varphi) \circ f = v^{n+2\over{n-2}} {\mathcal L}_g(v\varphi \circ f)$$
($\varphi = 1$ \ and \ $f = Id$ \ gives the Yamabe equation)}.\\

Consequently, we may deduce from the above, if \ $v$ \ is any positive solution of
\begin{equation}
{\mathcal L}_{g_0} v - {n-2\over{4(n-1)}} K v^{n+2\over{n-2}} = 0,
\end{equation} 
then \ $K$ \ is precisely the scalar curvature of the metric \ $g = v^{4\over{n-2}} g_0$.\\ Moreover, condition \ 
$g$ \  complete requires that \ $lim_{x\rightarrow \Lambda} v(x) = \infty$.\\

  The interest of the problem was
underlined by Schoen; he first proved the existence of weak solutions of Yamabe problem on \
$(S^n,g_0),$ \  which have a prescribed singular set \ $\Lambda_k$ \ .\\ In particular, this latter
point allows us to prove the existence of metrics \ $g$ \ conformally equivalent to \ $g_0$ \ , 
which have positive constant scalar curvature and are complete on \ $S^n - \Lambda_k$ \ . 
We refer to [M-Pa] for more details about the general problem, and for a
simpler version of the Schoen result.\\

The set \ $\Lambda_k$ \ may be replaced by \ $F(\Lambda_k) = \Lambda'_k$ \ , for
any conformal transformation. In particular, we may suppose that 
the points \ $p_i$ \  sum to zero, considered as vectors in the euclidian space
and then \ ${\Lambda}_k$ \ is contained in an equatorial subsphere \ $S^k
\subset S^n$.\  By using a reflection argument, one proves there is no complete (positive
regular) solution existing for the value \ $k = 1$ .\\
 One remarks a similarity between this singular Yamabe problem and
the theory of embedded, complete, constant mean curvature surfaces with \ $k$ \
ends in \ $I\!R^3$\ .\\
\ For the case \ $k =
2$\ , a family of such surfaces was discovered by C.
Delaunay. For this reason, it was called "Delaunay metrics" the analogous family of solutions for the
corresponding singular Yamabe problem. Some authors also called them  "Fowler metrics" in relation with radial solutions of Emden-Fowler equations.\\

 It is known the moduli space \ $ \mathcal{M}_\Lambda $ \ , which is the set of all
smooth positive functions \ $u$\ to the problem

\begin{equation}
 \qquad 4{{n-1}\over{n-2}} \Delta_{g_0} u + R(g_0) u - R(g)
u^{{n+2}\over{n-2}}  = 0
\end{equation} 
 $$g=u^{4\over {n-2}}g_0 \quad \mbox{is complete on} \quad
S^n-\Lambda\ ,
\qquad R(g) = \mbox{constant} > 0.$$ 

 This space is, in general a \ $k$\ -dimensional real
analytic manifold, see [M-P-U]. 
In the general case where \ $\Lambda = S^p$ \ (the p-subsphere where $1 \leq  p \leq  {n-2\over 2})$, \
 these solutions belong to an infinite dimensional space. F. Pacard [Pa] has shown the existence of
a positive complete solution in the case\ $ p = {n-2\over 2}. $ \ Nevertheless, the moduli space is
 only explicitly  determined if
\
$dim
\Lambda = 0,
\quad i.e.
\quad
\Lambda =
\Lambda_k$.\\  
In the case \ $k = 2$,\  we remark that the space
\quad $ \mathcal{M}_\Lambda $ \quad may be identified with any other, when \
$\Lambda$ \ contains just two elements. This space, denoted by \ $ \mathcal{M}_2$
\  which contains all the pseudo-cylindric metrics, may be identified with the open set \quad $\Omega \in
I\!R^2 $ \ . The frontier of \ $\Omega$ \  is the homoclinic curve \quad $\gamma_{b_0}$\ , which will be
described below.\\

 Since this problem is conformally invariant, we may choose the set \ $\Lambda_2
= (p,-p)$\ ; \ $-p$ \ is the antipodal point of \ $p$ .\\

Let \ $\Pi : S^n - \{pt\} \rightarrow IR^n$ \ the stereographic projection. We know that 
the converse \ ${\Pi}^{-1}$ \ defined by  $$\Pi^{-1} (x) = [ {2x\over{|x|^2 + 2}} , 
{ |x|^2 - 1\over{|x|^2 + 1}}]$$  is a conformal diffeomorphism, 
and we have \ $(\Pi^{-1})^\star (g_0) = {4\over{(|x|^2 + 1)^2}} dx^2$, \ where 
\ $ dx^2$ \ is the euclidean metric of \ $IR^n$.\\
Thus, the resolution of  (4) is equivalent to the resolution of
\begin{equation}
\Delta u + {n(n-2)\over 4} u^{n+2\over{n-2}} = 0
\end{equation}
on \ $IR^n - {\tilde\Lambda}_k,$ \ where \ ${\tilde\Lambda}_k = \Pi(\Lambda_k).$
 \ Moreover, \ $u$ \ have to verify\\ 
 \ $lim_{x \rightarrow {\tilde\Lambda}_k} u(x) = \infty$.
\subsection{Existence of metrics with two singular points } 
Examine now the case  \ $k = 2$.\ In this case, it will be possible to descrive all the 
Yamabe metrics.\ It is known that if \ $u$ \ is a solution on \ 
$S^n - \{p_1,p_2\}$, \ then it is invariant with respect to any conformal transformation
fixing these points. If these points are antipodal (as we may assume) 
then \ $u$ \ is rotationally invariant. It is then convenient to project the sphere 
\ $S^n$ \ to \ $IR^n$, \ from one of the singular point, \ $p_1$, \ so that \ $p_2$ \
sent to \ $0 \in IR^n$.\\ In this case, we have choosen \ ${\tilde\Lambda}_k = (0, \infty)$.\\
On the other hand, Cafarelli, Gidas and Spruck [C-G-S] proved that all solutions of (7)
on \ $ IR^n - \tilde{\Lambda}$ \  are radial. We then define the change
$$u(x) = |x|^{2-n\over 2} v(log{1\over{|x|}}).$$ In fact, this change corresponds 
to the conformal diffeomorphism \ $$\beta : IR^n - \{0\} \rightarrow IR\times S^{n-1}$$
where \ $ \beta(x) = (log{1\over{|x|}}, {x\over{|x|}}$).\\
Consequently, \ $v : IR \rightarrow IR^+$ \ is a solution of the autonomous ODE
\begin{equation}
{d^2\over{dt^2}} v  - {(n-2)^2\over 4} v + {n(n-2)\over 4} v^{n+2\over{n-2}} = 0
\end{equation}
The analysis of (8) show us, there is only one center say \ $(\alpha , 0)$ \ corresponding        
to the constant solution \ $\alpha = ({n-2\over n})^{n-2\over 2}$. \quad All the periodic orbits
 denoted by\ $\gamma_c(t)$, \  are parametrized
by \ $(v_c (t), {d\over{dt}}v_c (t)), \ c$ \ verifying \ $0 < c < c_0$.\ They are  surrounded by the homoclinic orbit \ $\gamma_{c_0}$.\ 
Here \ $v_{c_0} (t) = (\cosh t)^{n-2\over 2}$. \ The period is also depending on \ $c$ \ : \ $T = T(c)$.\\
The metrics $$g = {v_c}^{4\over{n-2}} (dt^2 + d{\xi}^2)$$ are called pseudo-cylindric type metrics. They are defined 
on the product manifold \ $(S^1
\times S^{n-1},dt^2+d{\xi}^2) $ \ a circle of length \ $T$ \  crossed with 
the standard
(n-1)-dimensional sphere. They belong to the conformal class of the cylindric metric
\quad $[dt^2+d\xi^2]$.

\subsection{Multiplicity of the pseudo-cylindric metrics}

R. Schoen  [Sc] gave a description of the family of metrics with constant scalar curvature in the conformal class of the product metric \ $g_0=dt^2+d\xi^2$\ on \ $S^1(T)\times S^{n-1}$\ 
the product Riemannian manifold  of a
circle of length \ $T$ \ and the euclidian (n-1)-sphere of radius 1.\\
Many other authors tried to descrive the conformal class  $[g_0]$, in order to complete the Schoen analysis concerning this problem, (see page 134 section 2 of [Sc]),
In particular, we refer to  [Ch1], and for a more geometrical interpretation we refer to [M-Pa].\\
The following result resume all the known properties of the pseudo-cylindric metrics, see [Sc], [Ch2], [M-Pa] and [Ch-W]

\bigskip

 {\bf Theorem (2-4)} \\
\quad {\it Let be \ $(S^1(T)\times S^{n-1},g_0=dt^2+d\xi^2)$ \ the product Riemannian manifold  of a
circle of length
\
$T$ \ and the euclidian (n-1)-sphere of radius 1. If the value of \ $T$\  satisfies the condition
\begin{equation}
T_{k-1} = {{2\pi (k-1)}\over{\sqrt{n-2}}} < T  \leq T_k = {{2\pi
k}\over{\sqrt{n-2}}},
\end{equation}
 where the integer \quad $k \geq  1$.\\
\ Then, there exist, in the conformal class \ $[g_0]$ \  k rotationally invariant pseudo-cylindric metrics\
$${g^j}_c = ({u^j}_c)^{4\over{n-2}}(dt^2+d{\xi}^2)\quad ,\quad j=0,1,2,..,k-1 \quad with \quad
{g^0}_c = g_0.$$ 
$k = 1$ \ only corresponds to the trivial product metric \ $C^{te} (dt^2+d\xi^2)$. }

\newpage

\section{
Remarks on the multiplicity of the Yamabe solutions }

\subsection{Yamabe metrics on manifolds}

Recall that for the standard sphere \ $(S^n,g_0)$\ , we only get a
trivial result: all Yamabe solutions are minimal for the functional \ $J(u)$\
and their sectional curvature turn to be constant. There actually exists a conformal
diffeomorphism between the manifolds \ $(S^n,u^{4\over{n-2}}g_0)$ \ and \
$(S^n,g_0)$ \ . In fact, all the metrics are obtained by pulling back
the canonical one, by a natural family of conformal transformations.\\ Among these metrics,
we have an explicit one parameter family of metrics (of constant scalar curvature
\ $R$\ )

\begin{equation}
 g_t = (\sqrt{1+t^2} + t\cos{\alpha x})^{-2} g_0\qquad \mbox{with  
} \qquad \alpha = \sqrt{R\over{n(n-1)}}.
\end{equation}
We can verify that all the functions \quad $u_t(x) = (\sqrt{1+t^2}
+ t\cos{\alpha x})^{-{{n-2}\over 2}}$\quad are solutions of the Yamabe equation (2) for \
$ R \equiv k$.\ These functions are minimizing for the functional, i.e.\ $J(u_t) = \mu (S^n,g_0).$\\

   Moreover, consider
 a Riemannian manifold $(M,g)$ \ with constant scalar curvature \ $R$ \ which in addition admits a conformal
isometric vector field \ $\mathcal X$.\\
 Then,  ${R\over{n-1}}$ \ belongs to the spectre of
the Laplace operator\ $\Delta_g$\ . By pulling back \ $g$ \ by the flow of \ $\mathcal X$ \, we get
a one parameter family of metrics which are conformal to \ $g$ \ and have the same scalar
curvature \ $R$ .

\subsection{The  degenerate case}
 
 $T= T_k$ \ corresponds to the degenerate case which has previously
been considered. The product metrics in general are not Morse critical points. Notice that
the Morse index of the solution \ $g = {u^{4\over{n-2}}} g_0$ \ is equal to the number of
eigenvalues of the Laplace operator
\ $\Delta_g$ \ in \ $[0,n]$ \ . The assumption that the solution is not degenerate implies that the
strict inequality \ $\lambda_1 > n$ \ holds. According to Lichnerowicz, we get \ $\lambda_1 \geq n$, \
whenever the metric \ $g$\ is Yamabe solution.\\  More generally, Schoen proves that, if all 
Yamabe solutions are non degenerate in a conformal class \ $[g_0]$ \ of a Riemannian compact manifold
\
$(M,g_0)$ \, then the functional has only a finite number of Morse critical
points. Thus, for a generic conformal class of metrics this hypothesis on \ $g_0$ \ is true, and \
$[g_0]$ \ contains at most a finite number of (non degenerate) Yamabe metrics. Moreover, the set of the
Yamabe metrics with finite energy is compact.

\subsection{Stability of the solutions and other multiplicity case}

 {\bf (i)} \quad Consider the stability problem of the pseudo-cylindric metrics\ ${g^j}_c =
({u^j}_c)^{4\over{n-2}}g_0$ ,\ we want to precise that among the \ $k$\ metrics given by the
above theorem, only the solution \
$u_T$
\ with fundamental period \ $T$\ is stable. Notice that the limiting period is the fundamental
period of the associated linearized equation of   $${d^2\over{dt^2}} u + (n-2) u.$$
The other
\ $(k-1)$\ solutions which have fundamental period \quad ${T\over j}$\quad with \quad $j =
2,...,k-1$ ,\quad are unstable.\\ By crossing the  point \
$(T_1,u_{T_1}(t))$ ,\ there is an ''exchange of stability'' between the two curves of solutions.
The trivial one becomes unstable after crossing this bifurcation point. Therefore, if we suppose \ $T > T_1$ \
the trivial solution \ $u \equiv \alpha$ \ does not minimize the Yamabe functional  $$ 
 J(u) = {{{4{{n-1}\over{n-2}}}\int_M\left|\nabla
u\right|^2 dv + \int_M R(x) u^2 dv}\over{\left(\int_M u^{2n\over{n-2}}dv\right
)^{n-2\over n}}}\ . $$ In fact, the infimum is attained by the
solution \ $u_T$ \ having the fundamental period \ $T$ \ .\\ If we get \ $T \leq  T_1$ ,\ the (unique)
solution \ $u \equiv \alpha$ \ minimizes the functional. To see that, it is enough to calculate
its value.\\ Indeed, let \ $V$ \ denotes the Riemannian volume of \quad $(S^1(T) \times
S^{n-1},g_0),$ 
\quad we get $$ J(u \equiv \alpha) = (n-1)(n-2)  \bigg(\int{dv}\bigg)^{2\over n}$$ where \
$$\alpha =  ({{n-2}\over{n}})^{{n-2}\over 4} \quad and \quad \bigg(\int{dv}\bigg) = T
vol(S^{n-1},d\xi^2) \ .$$  If we suppose $$T \leq  {2\pi\over{\sqrt{n-2}}}\ ,$$  obviously, we get 
$$J(u \equiv \alpha) < n (n-1) [vol(S^n)]^{2\over n}\ .$$  
Furthermore,  notice that only the non trivial solution 
\
$u_T$\ tends to  \quad $u_0(t) = (cosh t)^{-{{n-2}\over 2}}$,\quad when \ $T$ \ becomes
\ $\infty$\ . The minimum of the Yamabe functional \ $\mu_T$\ has to tend to \ $\mu(S^n,g_0) =
n(n-1){[vol(S^n)]}^{2/n}$, \ when \ $T$\ tends to \ $\infty$\ . \\

{\bf (ii) }\quad The multiplicity of the Yamabe problem may occur for other locally conformally
flat manifolds. Consider the real projective space: among the pseudo-cylindric metrics\
$g = u^{{4\over{n-2}}} g_0$,
\ we examine those  satisfying the condition  \quad $u(t,\xi) = u(-t,-\xi);$  \  any such solution
is independent of \
$\xi$\ . Thus we obtain  solutions to the Yamabe problem on the punctured real projective
space \ $\!R{P^n} - \{p\}$  \ by identifying any solution on this space with a solution on \ $S^n -
\{p,-p\}$,\ invariant under the antipodale reflection.

\subsection{The best constant of Sobolev and Aubin theorem}
There is another problem related on the multiplicity problem of Yamabe solutions:
 to determine the best constant of the Sobolev inclusions. Aubin [A1] has shown
that for a positive number \ $ \displaystyle{A \geq  K^{2}(n,2) = {4\omega_n^{-{2\over
n}}\over{n(n-2)}}}$ \  on the compact Riemannian manifold 
\ $(M,g)\quad of \ dimension\ n>2$, \ where \ $\omega_n$ \ denotes the volume of the standard
sphere, there exists a constant
\
$B$ \ such that, for  all \ $ u \in H^2_1(M,g)$ ,\ we have the following Sobolev inequality 
\begin{equation}
\Vert u\Vert^2_{L^{2n\over{n-2}}} \leq  A \Vert \nabla u\Vert^2_{L^2} + B \Vert
u\Vert^2_{L^2}\ .
\end{equation}
 The existence of the second constant  \ $B$\ seems
play a role in the Yamabe multiplicity problem and in the Nirenberg problem.
Hebey-Vaugon [H-V] have
 calculated the best little constant \ $B$\ denoted by \
$C(M,g)$.\ For this constant precisely, the inequality (11) is verified by any function \ $ u \in
H^2_1(M,g)$, \ on the quotient manifold of the standard sphere \
${S^n}/G$,\  where \ $G$\ is a cyclic group of isometries on the
sphere.\\ Moreover, one has estimated the best constant for a conformally flat manifolds  \ $S^1 \times S^{n-1}$\ 
endowed with a Riemannian metric product. More precisely, they find  $$C(S^1(T) \times S^{n-1}) \leq 
\displaystyle{{4\omega_n^{-{2\over n}}\over{n(n-2)}} T^{-2} + {n-2\over n}
\omega_n^{-{2\over n}}} $$
 where \ $T$\ is the circle radius.\\
Concerning the first constant \ $A$ \ . Notice that \ $A$ \ may attain the value
\ $K^{2}(n,2)$ \ if the manifold \ $(M,g)$ \ has constant sectional curvature.\\ According to
implicit function theorem, the Yamabe equation (4) is not locally invertible on \ $(S^1(T) \times
S^{n-1})$ \ near the trivial solution \ $u \equiv \alpha\ ,$ \ and $$\lambda = ({2\pi\over{T}})^2 =
n-2$$ is an eigenvalue, which yet verify the inequality $$\lambda < n
\bigg({\omega_n\over{\int{dv}}}\bigg)^{2\over{n}}\ .$$   
Obviously, the constant \ $K^{2}(n,2)$ \ is not attained, and there exists  \ $ u \in H^2_1(M,g)$ \
such that $$\Vert u\Vert^2_{L^{2n\over{n-2}}} > K^{2}(n,2)\bigg( \Vert \nabla u\Vert^2_{L^2} +
{(n-2)^2\over 4}
\Vert u\Vert^2_{L^2}\bigg)\ ,$$ 
since we have \quad $\mu(S^1(T) \times S^{n-1}) < \mu(S^n,g_0) = K^{-2}(n,2)$\\
Thus, we have proved the following.\\

{\bf Proposition (3-3)}\\
{\it Consider \ $V_n$\ a Riemannian compact manifold ($n \geq 3$), and the equation on \ $V_n$
$$ \Delta\phi + k\phi = k(1 + f) \phi^{n+2\over{n-2}}.$$
If, for \ $k = k_0$\ this equation is not locally invertible on \ $\phi_0 = 1$,\ then \ $\lambda =
{4{k_0}\over{n-2}}$\ is an eigenvalue of the Laplace operator.\\  
Moreover, if the eigenvalue verify\quad $\lambda \geq  n
{\displaystyle \bigg({\omega_n\over{\int{dv}}}\bigg)^{2\over{n}}}$,\quad then the Yamabe problem have
at least two solutions. Otherwise, \ $V_n$\ is conformal to the standard sphere.}\\

This proposition puts right and completes Theorem 11 of Aubin [A1], in including metrics like as the pseudo-cylindric metrics   (with respect to  notations of [A1])

\newpage

\section{
 Curvature properties}
\subsection{The Derdzinski metrics and the harmonic curvature}

\quad Let us consider a compact Riemannian manifold \ $(M,g)$\ , \ $\mathcal R$ \ is the 
curvature associated to the metric \ $g\ , \quad r = Ric(g)$ \quad is its Ricci tensor, and \ $D$\ is
the Riemannian
 connexion. We denote by \ $ \delta
\mathcal R$\ its formal divergence. This curvature is harmonic if \ $\delta {\mathcal R} = 0.$\ According to
the second Bianchi identity $$\delta {\mathcal R} = -d(Ric(g)) ;$$ we can express it in local
coordinates \
$D_k{\mathcal R}_{ij} = D_i{\mathcal R}_{kj}$.\quad In this case, \ $Ric(g)$ \ is a Codazzi tensor . In
particular, any Riemannian manifold with Ricci parallel tensor \ $Dr = 0$ \ (i.e. \ $D_i{\mathcal
R}_{kj} = 0$ ) \ has an harmonic curvature. Indeed, the Levi-Civita connexion \ $D$ \ is a
Yang-Mills connexion on the tangent bundle of \ $M$\ . In this way, the connexion \ $D$ \ is a
critical point of the Yang-Mills functional $$ {\mathcal Y} \mathcal{M} (\nabla) = {1\over 2} \int_M
\vert\vert R^\nabla\vert\vert dv,$$ where \ $R^\nabla$ \ is the curvature associated to the
connexion \
$\nabla$ \ . Notice that the Riemannian curvature must be harmonic for all Einstein manifolds and
for all conformally flat manifolds with constant scalar curvature; this
result can be deduced from the orthogonal decomposition of the curvature tensor. Moreover, the
condition of the curvature  harmonicity is, in a way , a generalisation of the Einstein
condition of the metric (\ $Ric(g)$\ is a constant multiple of \ $g$\ : \ $Ric(g) - {R\over{n}}g =
0$\ ). In particular, this fact shows that every Einstein metric must have a parallel Ricci tensor
\
$Dr = 0$.
\  In general, the latter property fails for a conformally flat metric with constant
scalar curvature, as we shall see below. However, in the 3-dimensional case, we get an identity
between harmonic curvature metrics and conformally flat metrics with constant scalar curvature.\\
 \ Bourguignon proved the converse for the compact manifold of dimension 4 and non
vanishing signature. But the signature assumption is necessary, as has been shown by 
Derdzinski\ [De].\ Moreover, Bourguignon has asked about the existence of compact metric with
harmonic curvature and non parallel Ricci curvature, as it was wellknown in the non compact case.\\
Derdzinski has given examples of metrics with harmonic curvature but non parallel Ricci tensor \
$Dr \neq 0$ \ , [De]. The corresponding manifolds are bundles with fibres \ $N$
\ over the circle \
$S^1$ \ (parametrized by arc length \ t\ ) with allowed warped metrics
\quad $dt^2 + h^{4/n}(t)g_0$ \quad on the product  \quad
$ S^1\times N $ .\quad Here, \ $(N,g_0)$ \ is an Einstein manifold of dimension \ $n \geq 3$, \  and the
function \ $h(t)$ \ on the prime factor is a periodic solution of the ODE, established by 
Derdzinski 
\begin{equation}
 h'' - {nR\over{4(n-1)}}h^{1-4/n} = -{n\over 4}C h \qquad \mbox{for some constant}
\qquad C > 0.
\end{equation}
 
 This function must be
non constant, otherwise the corresponding metric has a parallel Ricci tensor.\\ Notice that the
manifolds \ $S^1\times N $\ are not conformally flat, unless \ $(N,g_0)$ \ has constant
sectional curvature. Berger has shown that, under the hypothesis, of non negative sectional curvature
on a compact Riemannian manifold with harmonic curvature, the Ricci tensor is parallel .  
\\ J. Lafontaine proved that the connected sum
\
$ Y = \# X_i$ \ have the
same properties, where \ $X_i$ \ are manifolds with negative constant
scalar curvature. We may find many other examples of such manifolds in the chapter 16 of A. Besse [B].

\subsection{Curvature property of the pseudo-cylindric metrics}
 More generally, Derdzinski established a classification of the compact n-dimensional Riemannian
manifolds\quad $(M_n,g), \ n \geq 3$, \quad with harmonic curvature. If the Ricci tensor \ $Ric(g)$ \
is not parallel and has less than three distinct eigenvalues at each point, then \ $(M,g)$ \ is
covered isometrically by a manifold  $$(S^1(T) \times N,dt^2 + h^{4/n}(t)g_0),$$ where
the non constant positive periodic solutions \ $h$ \ verify the equation (11). Here \ $(N,g_0)$\ is a (n-1)-
dimensional Einstein manifold with positive (constant) scalar curvature.\\
 We get the following result, which gives another metric sharing this curvature property. See [Ch2]. \\

  {\bf
Theorem (4-3)} \\ {\it Consider the product manifold\quad $(S^1(T) \times S^{n-1},g_0)$\quad. Under
the condition  $$T(c) \geq T_1 = {{2\pi}\over{\sqrt{n-2}}},$$
  on the circle length, the Riemannian curvatures of
the associate pseudo-cylindric metrics\quad $g_c = {u_c}^{4\over{n-2}}g_0$ \quad are harmonic and their Ricci
tensors are non parallel.\\ Moreover, any pseudo-cylindric metric may be identified to a Derdzinski metric
up to a conformal transformation.}\\

\smallskip
Return now to the Derdzinski metrics descrived above. We remark that, in the locally conformally flat
case of this product, the factor \ $N$ \ may be identified with the standard sphere \ $S^{n-1}$\ .
Hence, these metrics must be conformally flat and will be conformal to the cylindric metric, as S.T.
Yau [Y] has shown \\ 

\smallskip

 {\bf Lemma (4-4) }\\
 {\it Any warped metric\quad $dt^2 + f^2(t)g_0$\quad on the product Riemannian manifold\quad $(S^1(T)
\times S^{n-1},dt^2 + g_0)$\quad must be conformally flat and conformal to a Riemannian
metric product \
$d\theta^2 + d\xi^2$ \ where \ $\theta $ \ is a $S^1$-parametrisation with length \quad $ 
\int_{S^1}{dt\over{f(t)}}$ }.\\

\smallskip
 Thus, any Derdzinski metric may be identified with a pseudo-cylindric metric up to conformal
diffeomorphism, let \ $F$ \ . Then the metrics are relied $$dt^2+{f^2}(t)d{\xi}^2 =
{F^*}({{u^j}_c}^{4\over{n-2}}(dt^2+d{\xi}^2)),$$ where \ ${{u^j}_c}$ \ are the pseudo-cylindric solutions
belonging to the (same) conformal class.\\
Indeed, for the metric \quad $dt^2 + f^2(t)d{\xi}^2$\quad  we can write $$ dt^2 + f^2(t)d{\xi}^2 = f^{2}(t)[({dt\over{f}})^2 +
d{\xi}^2].
$$After a change of variables and by using the conformal flatness of the product metric, we get $$dt^2 + f^2(t)g_0 = \phi^{2}(\theta)
[d\theta^2 + d{\xi}^2]$$ which is conformally flat. \\ This result is not surprising, because any manifold carying a warped metric
product
\quad
$(S^1 \times N,dt^2 + h^{4/n}(t)g_0)$\quad with harmonic curvature is not conformally flat
unless \ $(N,g_0)$\ is a space of constant curvature. This manifold must be locally conformally
equivalent to the trivial product \quad $S^1\times N.$ \\  This product will be conformally flat
only if
\ $N$\ has constant sectional curvature.
 
\subsection{
 Remark concerning the parallelism of the Ricci tensor}
 Parallelism property of the Ricci tensor have an interest particularly for the conformally flatness case.
Indeed, consider any Riemannian manifold \ $(M,g)$\ with a parallel Ricci tensor. This implies in particular, that its Weyl tensor is
harmonic:\quad
$\delta W = 0$\quad (in local coordonates $D_k W^{h}_{ijk} = 0).$\quad For a metric \ $\tilde g$\ in the conformal class and with
harmonic Weyl tensor \
$[g]$\ :
\
$\tilde g = e^{2\rho} g,$\ we get the following equality$$\tilde{\delta}{\tilde W} = \delta W - {1\over2}(n - 3) W
(\nabla\rho,.,.,.).$$ Since, the corresponding Ricci tensor is parallel, we then obtain in local coordonates 
 $$(n
- 3) W^{l}_{ijk} D_{l}\rho = 0.$$ Thus, when
\ 
$(M,g)$ \ is not conformally flat, then \quad
$\rho = 0$\quad necessarily\ (see "The Ricci calculus" of Schouten).

\newpage

 \section{
 On the existence of Derdzinski metrics}
\subsection{Analysis of the Derdzinski equation}

The above identification of the two family of metrics must reveal a similarity between the equations (8) and (12), since they have
been the starting point of our study. Indeed, by considering another time the ODE point of view, we are able to analyse all the
solutions of equation (12), like as we have made in the preceding section 2. Recall that they correspond to metrics with \quad $ Dr
\neq 0$.\   

So, this equation may have a non constant, positive periodic solution.\\ To that end, we
make a change of variables
\
$h(t) =
\alpha f(\beta t)\ ,$
\ where 
$$ \alpha =({(n-1)C\over{nR}})^{-n/4}\qquad {\mbox and} \qquad \beta = \sqrt{{nC\over 4}}.$$
Note that the constant \ $C$ \ must be positive to ensure that the equation (12) to have a periodic
non constant solution. This change gives the equation
\begin{equation}
  f'' - f^{1-4/n} + f = 0.
\end{equation}
We can verify that this equation satisfies the lemmas (2-6) and (2-7) given below in the last section. They will be used to
complete the analysis of the equation (8).\\ Indeed, we get the functions \\ $g(f) = f -
f^{1-{4\over n}} \quad;\quad g'(f) = 1 - (1-{4\over n})f^{-{4\over n}} \\ and \quad g''(f) =
{4\over n}(1-{4\over n})f^{-1-{4\over n}}$. \\ Notice that the hypothesis \quad $g''(f) \geq 0$\quad
is only satisfied in dimension \ $n \geq 4$.\\ Then, we calculate  $$\Delta(f) =
(f-1)[g''(1)g'(f) - g'(1)g''(f)].$$ \quad We get 
$$\Delta(f) = (f-1) {4\over n}(1-{4\over n}) [1 - f^{4\over n} + {4\over n}(f^{-{4\over n}} -
f^{-1-{4\over n}})].$$
It follows 
$$\Delta(f) = {4\over n}(1-{4\over n}) f(1-f^{-1})[1 - f^{4\over n} + {4\over n}f^{-{4\over n}}(1 -
f^{-1})],$$ 
which is positive.\\ Hence, we deduce the increase of the period function depending on the energy of
the equation (12). So, this fact allows us to determine the exact number of these Derdzinski
metrics. We first remark, that the bifurcation points of the solution family are \ $$(T_k,u_k) \quad
,\ where
\quad u_k = ({(n-1)C\over{nR}})^{-n/4}\quad and \quad T_k = {{2\pi}\over{\sqrt C}}.$$ 

\subsection{Existence conditions of the Derdzinski metrics}

 It appears from the above analysis that, a non constant, periodic solution of the equation (12)
exists only if the circle length satisfies the following condition
\begin{equation}
  T \geq {{2\pi}\over{\sqrt C}}.
\end{equation}
Notice that, we get an infinity of solutions . All are obtained by rotation-translation of the variable.\\ Moreover, in the case
where \ $T$ \ satisfies the double inequality 
\begin{equation}
2\pi {{(k-1)}\over{\sqrt C}} \leq T < 2\pi{k\over{\sqrt C}} \quad k \ {\it is \  an \  integer} \  > 1, 
\end{equation}
 then the
equation (8) may admit \ k \ rotationnally invariant distinct solutions. Hence, we have improved the lemma 1 of
Derdzinski [De] (see also Chapter 16 of A. Besse [Be]) \\

 {\bf Theorem (5-3) }\\
 {\it Consider the warped product \quad $(S^1(T) \times N,dt^2 + h^{4/n}(t)g_0)$\quad . 
The non constant positive periodic solution \ $h(t)$\ satisfies the equation (12) and \ $(N,g_0)$\
is a (n-1)-dimensional Einstein manifold with positive (constant) scalar curvature \ $R$\ . Then,
if the circle length \ $T$\ satisfies the following inequalities $$
2\pi{{(k-1)}\over{\sqrt C}} \leq T < 2\pi{k\over{\sqrt C}}\ , \qquad \mbox{where\  k \ is an integer} > 1,$$ there exist at least \ $k$\ rotationnally invariant warped metrics \\ $dt^2 + h^{4/n}(t)g_0$ \ on the product manifold \ $S^1(T)
\times N.$
\\ Moreover, these metrics have an harmonic curvature. Their Ricci tensor are non parallel only if \ $
 T \geq {{2\pi}\over{\sqrt C}}$\ .\\
 Conversely, if the warped metric of the type \ $dt^2 + h^{4/n}(t)g_0$ \ on the manifold \ $S^1(T)
\times N$
 has harmonic curvature, then the function \ $h(t)$ \ on the circle satisfies an ODE (12).}\\

\smallskip

The condition \ $
 T \geq {{2\pi}\over{\sqrt C}}$\ ,\ implies the existence of a non constant solution \ $h(t)$\ .
Otherwise, we have seen above that a trivial product have a parallel Ricci tensor.\\

 On the other hand, we remark that any warped metric  \ $dt^2 +
h^{4/n}(t)g_0$ \ defined by the lemma (4-4) on the product manifold \ $(S^1\times
S^{n-1},dt^2+d\xi^2) $\ is conformal to a Riemannian metric product $d\theta^2 + d\xi^2.$ \ Here \
$\theta $ \  is a new $S^1$-parametrisation with length \quad ${\displaystyle
\int_{S^1}{dt\over{h^{2/n}(t)}}}$ .\\
 Furthermore, we have seen in a previous paper [Ch1], there exists a analytic deformation
in the conformal class \ $[g_T]$\ , of any warped metric \ $g_T = dt^2 + f^2(t)d\xi^2$\ on the
manifold
\
$(S^1\times S^{n-1},dt^2+d\xi^2),\quad n=4 \ or \ 6$\ , and satisfying the length condition  $$
{\displaystyle T = \int_{S^1}{dt\over{f(t)}}}\ .$$ Notice that the product metric (under the
length condition) belongs to the conformal class\ $[g_T].$\  This analytic family of metrics depends
on two parameters \
$(g_{\alpha,\beta})$ \ . They all have a constant positive scalar curvature, and satisfy the
condition: \ $g_{\alpha,0}$ \ is the warped metric \ $g_T$\ .\\ Moreover, by using the same
argument of [Ch], we are able to extend the latter result to dimension 3. More precisely, a
metric on the euclidean space  \ $I\!R^3$ \ for which the rotation group \ $SO(3)$ \ acts by
isometries, is in fact a warped product as \ $g_T = dt^2 + f^2(t)d\xi^2$\ , where \ $d\xi^2$ \ is
the standard metric on the sphere \ $S^2$ \ . We can verify that its scalar curvature is $$R = {{2 -
2 {f'}^2 - 4 f f''}\over{f^2}}\ .$$ This fact and the connection with the dimension
\
$n$\ will be clarified below, in the next section.\\

Moreover, we know that every conformally flat manifold\ $(M,g)$ \ admits a Codazzi tensor which is
not a constant multiple of the metric. Let \ $b$ \ a such symmetric 2-tensor field; suppose \ $b$
\ has exactly 2 distinct eigenvalues\ $\lambda \ , \ \mu$ , has constant trace\ $tr_g(b) = c$ \ and
is non parallel. Following [De], these conditions give locally  
\begin{equation}
M = I\times N \quad
\mbox{with  allowed  metric}
\quad g = dt^2 + e^{2\psi} g_N
\end{equation}
 and 
\begin{equation}
 b = \lambda dt^2 + \mu e^{2\psi} g_N, \quad
\lambda ={c\over n} + (1-n) c e^{-n\psi}, \quad \mu = {c\over n} + c e^{-n\psi}\ .
\end{equation}

Conversely, for
any such data and for an arbitrary function \ $\psi$ \quad on \quad $I$ \ (16) defines a
Riemannian manifold \ $(M,g)$ \ with Codazzi tensor \ $b$ \ of this type. If we assume that the
positive function \ $h = e^{{n\over 2}\psi}$ \ satisfies   Equation (12)
$$ h'' - {nR\over{4(n-1)}}h^{1-4/n} = -{n\over 4}C h \qquad \mbox{for some 
constant}
\qquad C > 0,$$ then the
Ricci tensor is precisely a Codazzi tensor and, thus it is non parallel.\\ Therefore, \ $(M,g)$ \ is
isometrically covered by \quad  $(\!R \times N,dt^2 + h^{4/n}(t)g_N)$  \quad where \ $(N,g_N)$ \
is an Einstein manifold with positive scalar curvature.

\newpage

\section{
 Singular points of the pseudo-cylindric metrics: a global estimate }
\subsection{A bound for the pseudo-cylindric fonctions}

\ In this part, we examine another point of view coming from the explicit expression of
these metrics \quad ${g^j}_c = {{u^j}_c}^{4\over{n-2}}g_0$\quad. We can sometimes solve
explicitly the Yamabe equation. In particular, when it is reduced to the equation (8)
considered above. \\ Our approach has not only a calculation interest. However, the direct resolution is very instructive. In
particular, it gives the nature of the singularities of the solutions ,
which seems to depend on the dimension value \
$n.$\ \\  More precisely, for the dimensional cases \ $n = 3, 4$ \ and \ $6$\ , we
are able [Ch1] to resolve the equation (8) by using the elliptic functions. Notice that an
elliptic function is a double periodic meromorphic function, having the same number of poles and
zeros in a parallelogramme of periods. For these dimensions, the pseudo-cylindric functions
${u^j}_c$\ which have (necessarely) a real period \ $T$,\ turn out to be meromorphic, and have only
polar singularities. Moreover, \ $T$\ is precisely the period of the integral curve of Equation
(8), which is an elliptic curve (see Whittaker-Watson [W-W]). \\
 We recall that any solution must have a singular point; it becomes infinite on
approach of a point \
$p \in
\Lambda$\ . Furthermore, this solution (which is positive), must have two singular points in \ $0 \
and
\
\infty$\ , as has been shown by Cafarelli-Gidas-Spruck [C-G-S]. They have given an estimates near a singular point \ $p \in \Lambda$ ,\ that is for\ $
r$\ near \ $0,$ \ $u$\ satisfies\\
$$C^{-1}\ {r^{2-n\over 2}} \leq u(r) \leq C\ {r^{2-n\over 2}}$$
for some constant \ $C$ ,\ which depends only on\ $n$.\ In others words the metrics \  ${g^j}_c =
({u^j}_c)^{4\over{n-2}}g_0$\ must be cylindrically bounded near \ $p \in \Lambda$.\\
 R. Schoen have extended  this estimate to the whole manifold \ $S^n - \Lambda_k$,\ (Theorem (4-3) of
[P]). In the sense, these metrics globally are cylindrically bounded. His proof consider two important remarks. Firstly, by a
sequence of dilated coordonates \
$z_i = \lambda_i y$,\ the correspondant solutions \ $v_i(y)$\ are uniformly bounded on a large
ball in \ $ I\!R^n$.\ Secondly, the Harnack inequality ensures the existence of a positive
uniform bound for \ $v_i$\ in any compact set. By another method, in using properties of some
classical properties of elliptic functions, we find again the Schoen result for \ $k = 2$, and give
another estimate concerning the derivative. We obtain the following result, [Ch1].\\

\bigskip
 {\bf  Theorem (6-1) }\\
 \ { \it Consider the  pseudo-cylindric metrics \ ${g^j}_c = ({u^j}_c)^{4\over{n-2}}g_0$\ 
on the product Riemannian manifold \
$(S^1\times S^{n-1},dt^2+d\xi^2). $\ For the dimensional values \ $n = 3, 4 \ and \ 6$\ only,  the corresponding pseudo-cylindric fonction\
${u^j}_c$,\ may be expressed in terms of elliptic (meromorphic) functions.\\ For these dimensions , we get the estimates
$${\displaystyle {r^{n-2\over 2}} u(r) \leq C_n(T)\qquad and\qquad r^{n\over 2} {{du(r)}\over {dr}} \leq C'_n(T)},$$ where the
constants\ $C$\ and\
$C'$\ depend on the period \ $T$.\\ 
For the other dimensions,  the corresponding pseudo-cylindric fonction\
${u^j}_c$,\ may be expressed in terms of automorphic functions.}  
\smallskip

\subsection{
 General remarks }
 {\bf (i) }\quad 
 Consider a pseudo-cylindric metric\  ${g^j}_c =
({u^j}_c)^{4\over{n-2}}g_0 $  with  $ j=1,2,..,k-1$ \ . By changing the
variable \
$t$\ into a new variable \
$r(t)$\ , defined by  $${\displaystyle {dr\over{dt}}  = u^{2\over{n-2}} \quad and \quad r(0) =
0}\ ,$$ Mazzeo-Pollack-Uhlenbeck [M-P-U] have given a more simple expression of the period. Their
method, which is purely geometric, allows a quite explicit solution of the equation (8). Notice
that this new variable represents geodesic distance along the \quad
$\nabla t$\quad integral curves. The function \
$u(r)$\ is still periodic, and the period\ $R(c)$\ represents the length of the geodesic. For the
dimension \ $n = 4$,\ the solution can easily be  expressed as inverse sinus function of the
elliptic Jacobi function considered above. Then, the new period becomes \
$R(c)
\equiv
\pi$.\\

{\bf (ii) }\quad Now consider the linear differential equation\ 
$$ g\Delta f = Hess(f) - f Ric(g) \qquad (E) $$\ 
where \ $Hess(f)$\ denotes the second covariant derivative of \ $f\ , \Delta f =
trace[Hess(f)]$\ and \ $ Ric(g)$\ is the Ricci curvature tensor. Since this
equation is linear, $f = 0$\ is the trivial solution. This equation was
originally derived from the linearization of the scalar curvature equation. The
scalar curvature is considered as a non-linear map \ $\mathcal R$ \ of the space of \
$C^\infty $ \ riemannian metrics whose derivates are $L_2$-integrable to the
space of\ $C^\infty$\ functions on the manifold \ $M$ .\ This map is
differentiable and has derivative  \
$$d{\mathcal R}_g(\phi)  = {d\over{dt}}\big |_{t=0}{\mathcal R}(g+t\phi) = -\Delta tr
\phi + divdiv
\phi - \phi.Ric$$ and the adjoint operator \ $d{\mathcal R}_g^* \quad of \quad d{\mathcal R}_g$ \ with
respect to the inner product is
$$d{\mathcal R}_g^*(f) = -g\Delta f + Hess f - f Ric(g).$$
Thus \ $d{\mathcal R}_g^*$ \ is not surjective when \ $d{\mathcal R}_g^*$ \ has a
non-trivial kernel.\\

 The equation (E)  may have another (non-trivial) solution. Let \ $(M,g)$ \ be
the product manifold \ $(S^1 \times S^{n-1},dt^2+d{\xi}^2)$ \quad where \quad
$(S^{n-1},d{\xi}^2) $\ is the standard (n-1)-dimensional sphere with scalar
curvature \ ${\mathcal R}_0$\ and \ $(S^1,dt^2)$\ is the circle of length \ $T$ \
parametrized by \ $t.$ \  Under the additional geometric condition
\begin{equation} 
T = 2\pi\sqrt{{n-1}\over{R_0}}\ ,
\end{equation} 

 the equation have a non trivial solution on \ $(M,g)$\ .\\ This
condition corresponds precisely to the degenerate case of the product metric, in the sense that the
Yamabe solution is a degenerate critical point of the functional \
$J(u)$ \ . Moreover, this example (with other non trivial ones) was used to give a counter-example
to the Fisher-Marsden conjecture. \\ Notice that this manifold is conformally flat and has a
parallel Ricci tensor.\\
\quad Proposition 6.5 in [M-P-U]  asserts that, among the pseudo-cylindric
metrics, only the cylindrical (trivial product) one gives a non trivial solution to Equation
(E)
$$(E) \qquad  g\Delta f = Hess(f) - f Ric(g)  $$
but under the additional length condition\quad $T = T_k = {{2\pi k}\over{\sqrt{n-2}}}$. 
 
\newpage
\section{
 Curvature of the asymptotic pseudo-cylindric metrics} 
 An interesting and natural question is to ask if any
complete Yamabe metric on \quad $S^n -\Lambda_k$\quad where\quad $\Lambda_k = \{p_1,p_2,...,p_k\}\quad k > 2$\quad shares the same
curvature properties of the pseudo-cylindric metrics. This problem is not easy, since we do not have
 enough indications on these metrics, except their asymptotic behaviour on approaching a point \ $p_l
\in
\Lambda_k$. 
\subsection{Characterization of metrics with k singularities}

To be more clear, we use a result of [C-G-S] concerning the estimate of 
positive solutions \ $u(t,\xi)$ ,\  on the cylinder of the following
equation
\begin{equation} 
 4{{n-1}\over{n-2}} \Delta_{g_0} u + R(g_0) u - R(g)
u^{{n+2}\over{n-2}}  = 0, 
\end{equation} 
 
 $$g=u^{4\over {n-2}}g_0 \quad complete \ on \quad
S^n-\Lambda_k
\quad and \quad R(g) = constant > 0\ .$$ 
 By using a reflection argument, they actually assert that [C-G-S], any
solution \ $u(x)$ \ of equation (19), with a singularity at \ $p$ ,\ is asymptotic to one of the
pseudo-cylindric functions near the point \ $p$ \ . Following [P], we deduce that any such solution \
$u(x)$ \ corresponds to a solution \ $u(t,\xi)$ \ of equation (8) on a domain of \ $I\!R \times
S^{n-1}\ .$\\ Then, the corresponding pseudo-cylindric metric is unique.

 Let \ $\bar g = [u(t,\xi)]^{4\over{n-2}} g_0$ \ be a Yamabe metric conformal to the 
standard metric and \ ${u^{j,l}}_c(t)$ \ be the pseudo-cylindric solution at a point \ $p_l \in \Lambda_k.$\ The index\ $j$\ is
corresponding to one solution having period \ $T$.\  Then, we obtain the following as a corollary of a result of [C-G-S]\\

 {\it {\bf Lemma (7-2)} \\ As \ $t \rightarrow \infty$ \ we get the following estimate of the
Yamabe solution on\quad $S^n -\Lambda_k$\quad  near the point\ $p_l$
$$ u(t,\xi) = {u^{j,l}}_c(t) + {u^{j,l}}_c(t) {\mathcal O}(e^{-\beta t})\ ,$$
where\ ${u^{j,l}}_c(t)$\ is the corresponding pseudo-cylindric solution and \ $\beta$ \ is a positive constant}\\

 In the other words, this solution \quad $u(t,\xi)$ ,\quad which is asymptotic to a
pseudo-cylindric function \quad ${u^j}_c(t),$ \quad around a singular point\ $p_l$,\  can be written as \ $t \rightarrow \infty$ \\
\begin{equation} 
  u(t,\xi) = {u^j}_c(t) + e^{-\alpha t} v(t,\xi),
\end{equation} 
  where \ $\alpha$\ is a positive constant, only depending on the cylindric bound and \ $v(t,\xi)$ \ is bounded. This fact was
originally pointed out by  [C-G-S]. But [P], [M-P-U] have made it more clear,
by giving a precise estimate of the decay rate \ $\alpha$ . It may be calculated in terms of
spectral data of the operator \ ${\Delta_\epsilon} + n.$ \\ Indeed, approximate solutions for (19) which are singular at the points \
$p_l \in
\Lambda_k$
\ may be constructed by superposing translated and dilated of the singular radial solutions. \ The \ $k$\ free parameters in the exact solution correspond to the \
$k$\ different dilational parameters. This
fact has been precisely generalized by [M-Pa] in a more general context, to
approximate solution not only with isolated singularities, but also when \ $\Lambda$ \ is a
k-dimensional submanifold ($0 < k \leq {{n-2}\over 2}$) .\\  
 So, the corresponding metrics on \quad
 $S^n -
\Lambda_k$\quad are complete and asymptotic at each singular point to a (unique) pseudo-cylindric metric.
Moreover, these ``asymptotic metrics'' are uniquely determined by the Dilational Pohozaev
invariant \ ${\mathcal D}(g,p)$ \ on any isolated singular point. [P] has calculated its exact
value \\
\begin{equation} 
 {\mathcal D}(g,p) = 4{n-1\over{n-2}}\ \omega_{n-1}\ [E(b) + {1\over n}],
 \end{equation} 
where \ $E(b)$ \ is the Hamiltonian energy etablished in (2-6). Thus, the Dilational Pohozaev
invariant\ ${\mathcal D}(g,p)$ \ determines the Hamiltonian energies of the model of the
translated pseudo-cylindric solutions for the metric \ $g$ \ at each point \ $p_l$.

\subsection{ Ricci tensor of the asymptotic pseudo-cylindric metrics}

We get the following
result\\

{\it {\bf Theorem (7-3)}\\
 Any metric \ $\bar g$ \ of positive constant scalar curvature which is conformal to the standard
metric \ $g_0$ \ on the manifold \ $S^n-\Lambda_k$ ,\quad where \quad $\Lambda_k =
\{p_1,p_2,...,p_k\}\quad k > 2$,\quad  has a non parallel Ricci tensor, except for the standard
one.}\\

\bigskip

{\bf Remark (7-4)}\\ Near a singular point \ $p_l, \ l=1,2,...,k,$\ any Yamabe solutions on the manifold \ $S^n - \Lambda_k$\  have
another property related on the period
\
$T_l$\ of the corresponding pseudo-cylindric solutions. Consider its asymptotic expression, as \ $t \rightarrow \infty,$  $$ \qquad  u(t,\xi)
= {u^{j,l}   }_c(t) + e^{-\alpha t} v(t,\xi).$$ Let \ $T_l$ \ the (minimal) period of the corresponding pseudo-cylindric function \
${u^{j,l}}_c(t)$. This function verify the equality $$ u(t+T_l,\xi) - u(t,\xi) = e^{-\alpha t} [e^{-\alpha T_l}v(t+T_l,\xi) -
v(t,\xi)].$$ Since \ $v(t,\xi)$\ is bounded, then for any \ $t > N$\ we have  $$|u(t+T_l,\xi) - u(t,\xi)| < \epsilon .$$

\newpage

\section{Proofs}

\subsection{Proof of Theorem (2-4)}
 We first remark there is a conformal diffeomorphism between the manifolds \quad
$I\!R^n -
\{0\}$
\quad and \quad $I\!R\times S^{n-1}$ \quad given by sending the point \ $x$ \ to  \ $(log\vert x
\vert, {x\over{\vert x \vert}})$. \  By using the stereographic projection, we see easily that the
manifold
\
$I\!R
\times S^{n-1}$ \ (which is the universal covering space of \ $S^1\times S^{n-1}$) \  is 
conformally equivalent to \ $S^n -\{0,\infty\}$ \ .\\ Notice that the manifold \ $S^n - (p,-p)$ \
with the standard induced metric, can be considered as the warped product \ $$]0,\pi[ \times
S^{n-1}\ ,
\quad\mbox{with allowed metric} \quad dt^2 + \sin^{2}t d\xi^2$$ It has been shown by
Cafarelli-Gidas-Spruck [C-G-S] ( using an Alexandrov reflection argument), that any solution of (2-1)
is in fact a spherically symmetric radial function (depending on geodesic distance from either \
$p$ \ or \
$-p$ \ ). Any solution of Equation (6) which gives a complete metric on the cylinder \ $\Re \times
S^{n-1}$ \ is of the form \ $u(t,\xi) = u(t)$ ,\ where \ $t \in I\!R $ \ and \ $\xi \in S^{n-1}$.\\
The background metric on the cylinder is the product \ $g_0 = dt^2 + d\xi^2 .$ \  For the
convenience, we assume the sphere radius equal to 1.\\ Therefore, the partial differential equation is
reduced to an ODE.\\ 
 The cylinder has scalar curvature \ $R(g_0) = (n-1)(n-2)$ \ and \ $R(u^{4\over{n-2}}g_0) =
n(n-1).$  Thus \ $u = u(t)$ \ satisfies  
\begin{equation}
 {d^2\over{dt^2}}u - {(n-2)^2\over 4}u + {n(n-2)\over 4}u^{{n+2}\over{n-2}} = 0.
\end{equation}
 It follows that a pseudo-cylindric metric (constant scalar curvature metric) on the
product \ $(S^1(T)\times S^{n-1} , g_0)$ \ corresponds to a T-periodic positive solution of (8) and
conversely. The analysis of this equation shows us, that it has only one center  \ $(\alpha,0)$ \ .
 This corresponding to the (trivial) constant solution  $$ \alpha = ({{n-2}\over{n}})^{{n-2}\over
4},$$

\ and all periodic orbits \ $\gamma_b(t)$ \ are surrounded by the homoclinic orbit \
$\gamma_{b_0}$ .\ This orbit is parametrized by \ $$(u_0(t),v_0(t))
\quad  \mbox{where}
\quad u_0(t) = (cosh t)^{-{{n-2}\over 2}}.$$
 We denote the coordinates of \ $\gamma_b(t)$ \ by \ $(u_b(t),v_b(t).$ \
 When the value \ $b$ \ satisfies the condition \ $1 < b < b_0$ \ , the
correspondant orbit is periodic ($b_0$\ corresponds to a periodic solution of null energy, as we
shall see below the definition) .\\ One may easily remark that two positive T-periodic solutions of
(8) having the same energy are translated, and thus give rise to equivalent
metrics on \
$(S^1(T)\times S^{n-1},g_0)$. \\  Note that the metric corresponding to the conformal factor
\quad
$u_0\ : \quad  \ g = {u_0}^{4\over{n-2}}g_0$ \quad is incomplete. Then, it is not a pseudo-cylindric
metric; for this reason the constant \ $b$ \ cannot attain the critical value \ $b_0$ \ .\\ To
prove Proposition (2-4), we need the following result\\

 {\bf Lemma (2-5) }\\ {\it The solution family \ $(T,u_T(t))$ \ of the ODE \ (2-4),
(where \ $T$ \ is the minimal period) has  bifurcation points on the values
\quad $(T,u_T(t))$ \ where \ $ T_k={{2\pi k}\over{\sqrt{n-2}}}$\ and \ $u_{T_k}\equiv\alpha$.  
In this family, there is a curve of non trivial solutions which bifurcates to the right of the
trivial one.}\\

\smallskip
 This lemma is a classical result of global bifurcation theory (Crandall-Rabinowitz [C-R]). Indeed,
let us consider a positive T-periodic solution: if \ $T\neq T_k$\ , then the linearized associate
equation is non-singular.\\ We may also
deduce from bifurcation theorem, applied to the simple eigenvalues problem, that there is an
unique curve of non trivial solutions near the point \ $(T_k,\alpha).$ \  In fact, this uniqueness is
global. The trivial curve is \ $u_T \equiv\alpha$ .\\ Moreover, \
${\displaystyle ({du\over{dT}})_{T=T_k}}$ \ is an eigenvalue of the linearized associate equation.
According to the global bifurcation theory, we assert that the non trivial curves turn off on the
right of the singular solution \ $(T_k,u_{T_k}(t))$ \ . Consequently, when \ $T$ \ varies, two non
trivial curves never cross.\\

 On the other hand, we find it more convenient to analyze the equation (8), by making a change of
variables \ $u(t) = \alpha w(\beta t)$ \ ,\ where\\
$$ \alpha = ({{n-2}\over{n}})^{{n-2}\over 4} \quad  \mbox{and} \quad \beta =
{{n-2}\over 2} .$$
This change transforms (21) into the equation
\begin{equation}
 {d^2\over{dt^2}}w - w + w^{{n+2}\over{n-2}} = w'' + g(w) = 0.
\end{equation}
If \ $w(t)$ \ is a solution of (22), then the Hamiltonian energy $$E(t) = {{w'^2(t)}\over 2} -
{{w^2(t)}\over 2} + (n-2){w^{{2n}\over {n-2}}(t)\over {2n}} = {w'^2\over 2} + G(w) - {1\over n}$$
is independent on \ $t$ \ . Note that the constant solution \ $w(t) = 1$ \ has energy \
 $E_1 = -{1\over n},$ \  which is the smallest value that \ $E$ \ can take. Denote by \ $w_b$ \ the
solution of (22) which has \ $w(0) = b > 1$\ and \ $w'(0) = 0$ \ . Then, \ $${\displaystyle E(b) =
-{{b^2}\over 2} + (n-2){b^{{2n}\over {n-2}}\over {2n}} = c - {1\over n}}$$ \ is the energy of this
solution . If 
\

$1 < b < b_0 \ where \ b_0 = {1\over \alpha}$ \ satisfies \ $E(b_0) = 0$ \ , then \ $w_b(t)$ \ is a
positive, periodic function of $t$ , oscillating around the value $1$ , we have \ $$G(w) = \int_0^w
g(u)du \geq  0
\quad \mbox{satisfying}
\quad G(a) = G(b) = c,$$
\  where
$a$ belongs to the interval \ $[0,1]$ \ . Hence, we get the expression of the period $$T(c) = \sqrt
2\int^b_a{{dw}\over{\sqrt{c - G(w)}}}$$\\

 Then, to complete the proof of Proposition (2-4), we also need the following result of Chow and Wang.\\

{\bf Lemma (2-6) }\\
{\it The minimal period  \ $T(c)$ \ of a (periodic, positive) solution\ $u_c$\ of the equation (22)
is a monotone increasing function of its energy, when \ $c \in [0,{1\over n}]$}\\

\smallskip

 Chow-Wang [C-Wa] also have calculated the derivative of the period function, and have found the
expression of the derivative
\\
$$ T'(c) = {1\over c} \int_a^b {{g^2(w) - 2G(w)g'(w)}\over{g^2(w)\sqrt{c - G(w)}}}dw.$$
To prove the preceding Lemma, we need their following result (see corollary (2-5) in [C-Wa]) \\

\smallskip
  {\bf Lemma (2-7) }\\
 {\it Consider the  hypothesis on the smooth function \ $g \ \quad  g(1) = 0 \ ,\ g'(1) < 0$ \ ,
 we suppose in addition that
 $$H(x) = g^2(x) - 2G(x)g'(x) +{{g''(1)}\over{3g'^2(1)}} g^3(x) > 0,$$ 
for all \ $x \in [0,b_0],  $\ .\ Then, \qquad $T'(c) \geq
 0$ \qquad  for all} \ $c \in [0,{1\over n}]$ . \\

\smallskip
 Moreover, they show that under the conditions,
$$ g''(x) \geq 0 \quad and \quad \Delta(x) = (x-1)\bigg[g'(x)g''(1) - g'(1)g''(x)\bigg] \geq  0
$$
 the function \ $H(x)$ \ satisfies the previous Lemma (2-7) (that is the corollary (3-1) in [C-Wa]).

 Now, we apply these results to the function \ $$g(x) = -x + x^{{n+2}\over{n-2}}.$$ \ So we
get obviously
\
$$g'(1) = {4\over{n-2}}\quad , \qquad g''(1) = {4(n+2)\over{(n-2)^2}},$$ \ then we calculate
$$\Delta(x) = (x-1)\bigg[({{(n+2)}\over{(n-2)}}x^{4\over{n-2}} - 1) -
{{4(n+2)}\over{(n-2)^2}}({4\over{n-2}}x^{{6-n}\over{n-2}})\bigg],$$
and we find
 $$\Delta(x) = {{4(n+2)}\over{(n-2)^2}}(x-1)\bigg[(x^{4\over{n-2}} - 1) +
{4\over{n-2}}x^{{6-n}\over{n-2}}(x-1)\bigg],$$
which is positive.

\subsection{Proof of Theorem (4-3) }
To prove this Theorem, we consider a local chart in this product manifold, where the local
coordonates in this chart are 
\quad
$(x^0,x^1,....,x^{n-1}),$  \ and \ $x^0 = t$ \quad is the coordinate corresponding to the circle
factor\quad
$S^1$. \quad We can write any metric tensor in the conformal class \ $[g_0]$\quad : \quad
$\bar{g} = {u}^{4\over{n-2}} g_0,$
\quad where the\ $C^\infty$\ function\ $u$\ is
  defined on the circle. Then we obtain the expressions in the local coordonates
system 
\begin{equation}
 \bar{g}_{00} = {u}^{4\over{n-2}} , \quad \bar{g}_{0i} = 0 , \quad \bar{g}_{ij} =
{u}^{4\over{n-2}} (g_0)_{ij}, \quad
\end{equation}
where
\quad $i , j = 1,2,...,n-1. $ \\ 
For the following, we refer to Kobayashi-Nomizu [K-N].\\
We get the Christoffel symbols\quad $\bar{\Gamma}^i_{jk}$ \quad associated to the metric \
$\bar{g}.$ \ They are defined by the expression $$ \bar{\Gamma}^i_{jk} = {1\over 2}
\bar{g}^{il} {{\partial \bar{g}_{jl}}\over{\partial x^i}} + {{\partial \bar{g}_{il}\over{\partial x^j}} -
{\partial \bar{g}_{jk}\over\partial x^l}}).$$
Then, we obtain  
$$\bar{\Gamma}^0_{jk} = - {2 \over{n-2}}{u' \over u} {\bar g}_{jk}, \quad \bar{\Gamma}^i_{j0} = {2
\over{n-2}}{u' \over u}{{\bar \delta}^i_j},\quad
\bar{\Gamma}^0_{00} = {2 \over{n-2}}{u' \over u}\ and \ \bar{\Gamma}^i_{00} = 0$$
The Ricci tensor components are 
$${\bar\mathcal R}_{ij} = \mathcal R_{ij} - {2\nabla_{ij} u \over{u}} + {2n \over{n-2}}{\nabla_i u \nabla_j
u \over{u^2}} - {2 \over{n-2}} {{\mid \nabla u\mid^2 + u\Delta u} \over{u^2}} g_{ij}.$$
We have  
 $${\bar{\mathcal R}}_{0i} = 0 , \qquad {\bar {\mathcal R}}_{00} = 2{{n-1}\over{n-2}} [ {u''\over u} +
{u'^2\over{u^2}} ].$$
The associated scalar curvatures are 
$$R(g_0) = (n-1)(n-2) \qquad and \qquad {\bar R}(u^{4\over{n-2}}g_0) =
n(n-1).$$ \quad So, the trace of the Ricci tensor gives the Yamabe equation
$${\bar g}^{ij}{\bar\mathcal R}_{ij} = {\bar{\mathcal R}} = u^{-4 \over{n-2}} (R + {4(n-1)
\over{(n-2)}}{\Delta u
\over{u}}),$$ and the function
\quad
$u $
\quad satisfies the above equation 
$$ u'' - {(n-2)^2\over 4}u + {n(n-2)\over 4}u^{{n+2}\over{n-2}} = 0.$$  

The covariant derivative of the Ricci tensor gives
$$D_k {\bar\mathcal R}_{ij} = {\partial{\bar\mathcal R}_{ij}\over{\partial{x^k}}} - \bar{\Gamma}^l_{ik}
{\bar\mathcal R}_{lj} -\bar{\Gamma}^l_{jk} {\bar\mathcal R}_{li},$$ 
where \ $D$ \ denotes the riemannian connexion associated to the metric \ $g$ .\\ 
 In particular we get  
$$D_0\bar{\mathcal R}_{00} = {d\bar{\mathcal R}_{00} \over {dt}} - 2\bar{\Gamma}^k_{00}\bar{\mathcal R}_{k0}
.$$ 
Thus, we obtain $$D_0\bar{\mathcal R}_{00} = 2{{n-1}\over{n-2}} {{u'''}\over u} +
2{{n+2}\over{n-2}}{u'u''\over{u^2}} - 4{{u'^3}\over{u^3}}.$$
By comparison with (2-4), \ $u$ \ should be a periodical function. So, we deduce that,
if the function
\
$u$
\ is non constant, then 
$$D_0\bar{\mathcal R}_{00} {\neq} 0.$$
The Ricci curvature of the pseudo-cylindric metrics is not parallel, except for the cylindric one. 

\subsection{Proof of Theorem (6-1) } 
\ Consider the orbits of the autonomous system associated to the equation (8)  \\ 
\begin{equation}
u'^2 = ({n-2\over 2})^2\bigg[u^2 - u^{2n\over{n-2}}\bigg] + \bar{c}
\end{equation}
these orbits are closed only if the energy constant satisfies the following condition
\begin{equation}
 -{2\over n}({n-2\over 2})^2({n-2\over n})^{n-2\over 2} < \bar{c} < 0.
\end{equation}
In connection with the energy constant \ $c$ ,\ used in the proposition (2-4), we obtain
the
  relation \ $$\bar{c} = -2\alpha{\beta}^2 c = -{(n-2)^2\over 2} ({n-2\over n})^{n-2\over
2} c,$$ \quad where the constants are as above\ $\alpha = ({n-2\over n})^{n-2\over 2}\quad
and\quad \beta = {n-2\over 2}$\quad (we recall that the constant \ $c$\ must varying in the interval
\ $[0,{1\over n}]$).\\
 We remark that for the n-values \ $n = 4\ or \ 6$ ,\ all the orbits are elliptic curves, which may
be parametrized only by 2-periodic meromorphic functions (i.e. elliptic functions).\\
 More exactly, a direct calculation gives the explicit solutions, for \ $n = 4$ \\
\begin{equation}
 u_{\bar{c}}(r) = {\displaystyle{1\over{\sqrt{2-\bar{c}^2}}} \ 
dn_{\bar{c}}({r+\beta\over{\sqrt{2-\bar{c}^2}}})} 
\end{equation}
 where \ $dn_{\bar{c}} $\ is the elliptic Jacobi
function\\
$$dn_{\bar{c}}(x) = 1 - {{\bar{c}}^2\over{2}} x^2 + {{\bar{c}^2}(4+{\bar{c}^2})\over{6}} x^4 + ...$$
 of real period \\
$$T(\bar{c}) = 2{\sqrt{2-\bar{c}^2}}\int_0^{\pi/2}{d\theta\over{\sqrt{1 - \bar{c}^2\sin^2\theta}}}.$$
These Jacobi functions \ $dn_{\bar{c}}(r)$\  satisfy the differential equation\\
$${({d\over{dr}} dn)^2} = (1 - dn^2)(dn^2 + \bar{c}^2 - 1)$$
 with \ $dn_{\bar{c}}(0) = 1$\ , Whittaker-Watson [W-W].\\

 We obtain analogous results for the case \ $n = 6$ ,\ where the solution may be
expressed in terms of the Weierstrass function\ $\wp(r)$. We recall the equation satisfied by the
Weierstrass function $$ \wp'^2(r) = 4\wp^3(r) - g_2\wp(r) - g_3 ,$$ and its expression is $$\wp(r)
=  {1\over{r^2}} + \sum_{0<p} a_{p}(g_2,g_3) r^p.$$ In particular we deduce that, \
$r^2\wp(r)$ \ is an entire function in an elementary parallelogram of periods. Then, it is bounded.\\
 So, we get the following relation \\
\begin{equation}
 u_{\bar{c}}(r) = {1\over 3} - \wp(r+\beta;g_2,g_3),
\end{equation}
where \ $g_2 = {4\over 3}$\quad and \ $g_3$ \ which depends on the parameter \ $\bar{c}:\ (g_3 = -{8\over27}-\bar{c})$\  satisfies
the condition\quad
${g_3} > - ({2\over 3})^3$.\\ This latter inequality, corresponding to the (25), insures us that
the function \ $\wp(r)$\ has a real period. Furthermore, if \ $r = 0$\ is a singular point, then we can take the value \ $\beta =
0$.\\ It is known that this elliptic function has a double pole at the origin and two zeros (in a elementary parallelogram of
periods). Thus, its  derivative has 3 zeros and a pole of order 3, [W-W].\\

Moreover, it is known when\ $r$\ increases in the interval \ $]0,{T\over2}[$,\ then the Weierstrass function\ $\wp(r)$\ decreases
monotonically in \ $]+\infty,\wp({T\over2})[$.\ So, we obtain explicitly for \ $n = 6$\ $$ r^2 u(r) \leq T^2 u({T\over2}).$$ Also, for
the dimension \ $n = 3$\ or \ $n = 4$, \ we get an analog $$r^{{n-2}\over{2}} u(r) \leq T^{{n-2}\over{2}} u({T\over2}).$$ \\

 Examine now the general case: to study the integral curves of our equation (24) and their
algebraic nature; it seems that it will depend on the parity of the dimension. For
this, we consider separately two hypotheses concerning the parity.\\ 
 
\quad {\bf (i)  }\ Suppose that the dimension \ $n$\ takes an odd value  \ $n = 2p+1.$\  Using the
change of function $$u = v^{-{n-2\over 2}}\ ,$$ we get \quad $u'^2 = ({n-2\over 2})^2
v^{-n} v'^2$\quad and, after the replacement, Equation (24) becomes\\
$$({n-2\over 2})^2 v^{-n} v'^2 = ({n-2\over 2})^2 [ v^{-(n-2)} - v^{-n} ] + \bar{c}\ ,$$ and we
obtain the analogous equation for \ $v$ 
\begin{equation}
 v'^2 = v^2 - 1 + ({2\over{n-2}})^2 \bar{c}\ v^n\ .
\end{equation}
 By referring to the geometry of algebraic curves (Shafarewitch [Sh], for any odd value \ $n \neq
3$\ , the integral curve of this equation is known to be represented by a hyperelliptic curves. Whose
coordonates are connected by the algebraic relation \\
$$ F(x,y) =  y^2 - x^2 + 1 - ({2\over{n-2}})^2 \bar{c}\ x^n = 0\ .$$\\
  Under condition (26) on \ $\bar{c}$\ , the integer \ $p$\ is precisely the genus of such
curves if the polynomial \quad $P(x) = x^2 - 1 + ({2\over{n-2}})^2 \bar{c}\ x^n$\quad is without
multiple roots. Otherwise, the genus fails (which corresponds to the extremal values of \ $\bar{c}$\
).\\ Moreover, as we may easily remark, the curve is elliptic  only for the value
\
$n = 3$\ . For this particular case, we get the equation $$(a) \qquad v'^2 = v^2 - 1 + 4 \bar{c}\
v^3.$$\\ This equation has a elliptic functions as solutions. More exactly, under the condition
(26) on \
$\bar{c},$
  $$ -{1\over{6\sqrt 3 }} < \bar{c} < 0$$ we get the solution $$ v_{\bar{c}}(r) = 
{1\over{\bar{c}}} [ \wp (r+\beta;g_2,g_3) - {1\over {12}} ],$$
where \ $\wp$\ is the Weierstrass function relative to \ $g_2 = {1\over{12}}\ and \ g_3 = \bar{c}^2 -
{1\over{216}}$. \ This function has two simples zeros in an elementary parallelogram of periods, one
of them is real. We choice \ $\beta$\ such that\ $wp(\beta) = {1\over{12}}$.\ So, we get \ $v_{\bar{c}}(0) = 0$.\ We deduce that the
function
\quad
$u_{\bar{c}} = v_{\bar{c}}^{-{1\over 2}}$\quad has a pole of order \ $1\over 2$. \\Otherwise, for any odd value of the dimension\ $n
> 3$\ , the corresponding integral curve of (29) does not admit a parametric representation in terms of
meromorphic functions, because it has genus \ $p > 1$\ .\\ Equivalently, the algebraic relation \
$F(x,y) = 0$\ cannot be uniformized by meromorphic functions in the finite plane. In fact, this
curve may be uniformized   $$x = A(r)\quad ; \quad y = B(r),$$
where \ $A$\ and \ $B$\ are automorphic functions which are meromorphic but have only a
restricted domain of existence with a natural boundary.\\

\quad {\bf (ii) - }\ Now, consider the case where the dimension\ $n$\ takes an even value. This seems
to be harder than the latter case; in particular determining the genus of the integral
curve is less easy.\ We use a more appropriate change
$$u = v^{-{n-2\over 4}}\ ,$$ we get
\quad
$$u'^2 =  ({n-2\over 4})^2 v^{-{n+2\over 2}} v'^2$$\quad and the equation (24) becomes \\
$$({n-2\over 4})^2 v^{-{n+2\over 2}} v'^2 = ({n-2\over 2})^2 [ v^{-{n-2\over 2}} - v^{-{n\over 2}} ]
+ \bar{c},$$ then we deduce the analogous equation for \ $v$
\begin{equation}
  v'^2 = 4v^2 - 4v + ({4\over{n-2}})^2 \bar{c}\ v^{n+2\over 2}.
\end{equation}

 To know the algebraic nature of the integral curve of this equation, we proceed as before. Remark
that, for the values \ $n = 4\ or \ 6$ ,\ the associated integral curve for this equation is
 an elliptic one (here the genus is \ $p = 1$\ ) . However, these two cases have been
treated before.\\ Otherwise, for any even value \ $n \geq 8$\ , the integral curve of (29) is a
hyperelliptic one, if its second member is without multiple roots.\\ The genus \ $p$\ of this curve
depends on the parity of \ $n\over 2$\ ; it is  \quad $\biggl[{n\over 4}\biggl]$ ,\quad the
integer part of\
$n\over 4$\ .\\ It was proved in [Sh], that the case where the polynomial \ $P$\ has even degree,
reduces to that of odd degree.\\ Moreover, this transformation preserves the genus.
For this reason the curves corresponding to \ $n = 4\ and\ 6$\ have  genus \ $p = 1$,\  likewise for
those corresponding to \ $n = 8\ and\ 10$ ,\ which have  genus \ $p = 2$\ , and so on.\\ Generally
speaking, every integral curve of (28) and (29) admits \ $2p$ abelian periods. One of them is
real, and this one is interesting for our problem.\\
On the other hand, except for the cases\ $n = 3,4$\ and\ $6$\ , the general solution of the
equations (29) and (30) may be expressed in terms of automorphic functions. These functions may be
regarded as a parametric representation of the corresponding integral curves.\\ It is known that this
solution is an automorphic function of a Fuchsian group of the first kind, all of whose
homographic substitutions are hyperbolic. The fundamental region (for the same genus \ $p$) is
contained in the interior of the limit circle.\\
Notice that every limit point is necessary an essential singularity for any (non degenerate)
solution.\\ In general, it seems to be very difficult to obtain explicit formulas, the theory of
elliptic functions is rather exceptional.
  
\subsection{Proof of Theorem (7-3)}
To prove this, we need the asymptotic property of these metrics. In particular,
the estimate (20) play an essential role.\\ Notice that, this  manifold \ $(S^n-\Lambda_k,\bar
g)$ \ is locally conformally flat. The condition of constant scalar curvature implies that its
Riemannian curvature is harmonic. We get the converse only for the dimension \ $n = 3$ .

\bigskip
We proceed as in the proof of Theorem (4-3 ). Consider a local chart \ $(t,\xi)$\ of a point in \
$S^n-\Lambda_k$ \ , here \ $t = x^0$\ and \ $\xi = (x^1,x^2,...,x^{n-1})$ .\\
 In this local
coordonates system, we obtain the components of the tensor metric $${\bar g}_{00} =
u^{4\over{n-2}}\quad , \qquad {\bar g}_{0i} = 0 \quad , \qquad {\bar g}_{ij} = u^{4\over{n-2}}
(g_0)_{ij}.$$

The Christoffel symbols \ ${\bar \Gamma} ^{i}_{jk}$ \ associated to the metric tensor\
${\bar g}$ are
$$ \bar{\Gamma}^{i}_{j0} = {2\over{n-2}}{1\over u} [{\partial u \over{\partial x^j}} +
{\partial u\over{\partial x^i}} \delta^i_j] \qquad {\bar \Gamma}^{0}_{00} = {2\over{n-2}}{1\over
u} {{\partial u}\over{\partial t}}$$

Then we get  components of the Ricci tensor $${\bar{\mathcal R}}_{0i} = - {2\nabla_{i0} u \over{u}} + {2n \over{n-2}}{1\over{u^2}}{\partial
u\over{\partial x_i}}  {\partial
u\over{\partial t}} ,$$ 
$$ {\bar {\mathcal R}}_{00} = 2{{n-1}\over{n-2}} [ {1\over u}{\partial
^2}u\over{\partial {t^2}} + {1\over{u^2}} ({{\partial u}\over{\partial t}})^2 ].$$

 \ $D$ \ denotes the Riemannian connexion associated to the metric \ ${\bar g}.$ \\ 
Considering now Lemmas (6-2) and (6-3) which imply that  $${{\partial u}\over{\partial x^i}} =
e^{-\alpha t}{\partial v\over{\partial x^i}} \quad and
 \quad \nabla_{i0} u = -\alpha e^{- \alpha 
t}{{\partial v}\over{\partial x^i}}
+ e^{-\alpha t}{{\partial^2 v}\over{{\partial x^i}{\partial
t}}}.$$

Recall that \quad $v = v(t,\xi)$\quad is a bounded function, and \ $\alpha$ \ is a positive
constant.\\ 

Thus, we can write 
$${\bar{\mathcal R}}_{0i} = e^{-\alpha t} f(t,x_1,x_2,...,x_{n-1})$$
where \ $f$ \ is bounded .\\ 

Consequently, we may calculate the components of the covariant
derivative
 $$D_0\bar{\mathcal R}_{00} = 2{{n-1}\over{n-2}} {{{u_c}'''}\over u_c} +
2{{n+2}\over{n-2}}{{u_c}'u''\over{{u_c}^2}} - 4{{{u_c}'^3}\over{{u_c}^3}} - e^{-\alpha t}
g(t,x_1,x_2,...,x_{n-1})$$
where \ $g$ \ is a bounded function.\\ 

Therefore, since \ $t$ \ is large, we conclude that 

$$D_0\bar{\mathcal R}_{00} {\neq} 0.$$

\newpage

{\it {\bf REFERENCES}}\\

[A1]  T.Aubin \ {\it Meilleures constantes dans le théorème d'inclusion de Sobolev }\quad 
J. Of Funct. Anal., 32, 148-174, (1979).\\ 

[Be]  A.L.Besse \ {\it Einstein Manifolds}\quad Ergebnisse der
Mathematik und ihrer Grenzgebiete, 3. Folge, Band 10. Springer. New-York, Berlin (1987).\\

[C-G-S]  L.Cafarelli, B.Gidas and J.Spruck \ {\it Asymptotic symmetry and local behavior of
semilinear elliptic equations with critical Sobolev Growth}\quad Comm. Pure and Appl. Math., vol
XLII (1989), 271-297.\\

[Ch-W] R.Chouikha and F.Weissler \ {\it Monotonicity properties of the period function and the number of constant positive scalar curvature metrics on $S^1\times S^{n-1}$}\quad
Prepubl.  Math., Universite Paris-Nord, n° 94-14, (1994).\\

[Ch1]  R.Chouikha \ {\it Famille de m{\'e}triques conform{\`e}ment plates et
r{\'e}gularit{\'e}}\quad C.R. Math. Reports. Acad Sc. of Canada, vol XIV (1992) 257-262.\\

[Ch2] A.R. Chouikha \ {\it On the properties of certain metrics with constant scalar curvature}\quad
 Diff. Geom. and Appl., Sat. Conf. of ICM in Berlin, (1999), Brno, 39-46.\\ 
 
[Ch3] A.R. Chouikha \ {\it Nonparallelism of the Ricci tensor of pseudo-cylindric metrics}\quad
Geom. Dedicata 93 (2002), 107--112.\\ 
 
[Ch4]  A.R. Chouikha \ {\it Existence of metrics with harmonic curvature and non parallel Ricci tensor}\quad 
 Balkan J. Geom. Appl. 8 (2003), no. 2, 21--30.\\

[C-Wa]  S.N.Chow and D.Wang \ {\it On the monotonicity of the period function of some second order
equations}\quad Casopi pro pestovani matematiky, 111 (1986), 14-25.\\

[C-R]  M.Crandall and P.Rabinowitz \ {\it Bifurcation from simple eigenvalues }\quad J. Functional
Analysis 8, (1971), 321-340.\\

[De]  A.Derdzinski \ {\it Classification of certain compact Riemannian manifolds with harmonic
curvature and non parallel Ricci tensor} \quad Math. Z., 172, (1980), 277-280.\\

[De1]  A. Derdzinski \ {\it On compact Riemannian manifolds with harmonic curvature } \quad Math.
Ann., vol 259, (1982), 145-152.\\

[H-V]  E.Hebey and M.Vaugon \ {\it Meilleures constantes dans le th{\'e}o{\`e}me d'inclusion de
Sobolev et multiplicit{\'e}s pour les probl{\`e}mes de Nirenberg et de Yamabe} \quad Indiana Univ.
Math. J., 41, (1992), 377-407.\\

[K-N]  S.Kobayashi and K.Nomizu \ {\it Foundations of Differential Geometry} \quad vol. 2. New
York-London, Interscience (1966).\\

[L-P] J.M.Lee and T.M.Parker \ {\it The Yamabe problem}\quad Bull. Amer. Math. Soc., 17, (1986)
p.37-91.\\

[M-P-U]  R.Mazzeo, D.Pollack and N.K.Uhlenbeck \ {\it Moduli spaces of singular
Yamabe}\quad J. Amer. Math. Soc. 9 (1996), no. 2, 303--344.\\

[M-Pa]  R.Mazzeo and F.Pacard \ {\it Constant scalar curvature metrics with isolated singularities}\quad MSRI, dg-ga server,  Duke Math. J. 99 (1999), no. 3, 353--418.\\

[Pa]  F.Pacard \ {\it The Yamabe problem on subdomains of even dimensional spheres}\quad
Topol. Methods Nonlinear Anal. 6 (1995), no. 1, 137--150.\\

[P]  D.Pollack \ {\it Compactness results for complete metrics of constant positive scalar
curvature on subdomains of $S^n$} \quad Indiana Univ. Math. J., vol 42, (1993), 1441-1456.\\

[Sc]  R.Schoen \ {\it Variational theory for the total scalar curvature functional for
Riemannian metrics and related topics} \quad Lectures Notes in Math., 1365, Springer-Verlag
(1989), 120-154.\\

[{Sc1}]  R.Schoen \ {\it On the number of constant scalar curvature metrics in a conformal class}
 \quad Differential Geometry: A symposium in honor of M. Do Carmo. Wiley, (1991), 311-320.\\

[Sh]  I.C.Shafarewitch \ {\it Basic algebraic geometry}\quad Springer-Verlag, (1977).\\

[W-W]  E.T.Whittaker and G.N.Watson \ {\it A course of modern analysis}\quad Cambridge, (1963).\\

[Y]  S.T.Yau \ {\it Remarks on conformal transformations}\quad J. Diff. Geom., 8, (1973), 369-381.

\end{document}